\documentclass[twosided,reqno]{amsart}

\usepackage{amsmath}
\usepackage{amssymb}
\usepackage{amsthm}

\theoremstyle{theorem}
\newtheorem{theorem}{\scshape Theorem }[section]
\newtheorem{lemma}[theorem]{\scshape Lemma}

\theoremstyle{definition}

\numberwithin{equation}{section}

\begin{document}

\title{Some identities of $q$-Euler polynomials arising from $q$-umbral calculus}

\author{Dae San Kim}
\address{Department of Mathematics, Sogang University, Seoul 121-742, Republic of Korea.}
\email{dskim@sogang.ac.kr}

\author{Taekyun Kim}
\address{Department of Mathematics, Kwangwon University, Seoul 139-701, Republic of Korea}
\email{tkkim@kw.ac.kr}

\maketitle

\begin{abstract}
Recently, Araci-Acikgoz-Sen derived some interesting identities on weighted $q$-Euler polynomials and higher-order $q$-Euler polynomials from the applications of umbral calculus
(See [1]).  In this paper, we develop the new method of $q$-umbral calculus due to Roman and we study new $q$-extension of Euler numbers and polynomials which are derived from
$q$-umbral calculus.  Finally, we give some interesting identities on our $q$-Euler polynomials related to the $q$-Bernoulli numbers and polynomials of Hegazi and Mansour.
\end{abstract}

\section{Introduction}

Throughout this paper we will assume $q$ to be a fixed real number between $0$ and $1$.  We define the $q$-shifted factorials by
\begin{equation}\label{1}
(a:q)_0=1, (a:q)_n=\prod_{i=0}^{n-1}(1-aq^i), (a:q)_{\infty}=\prod_{i=0}^{\infty}(1-aq^i).
\end{equation}

If $x$ is a classical object, such as a complex number, its $q$-version is defined as $[x]_q=\frac{1-q^x}{1-q}$.
We now introduce the $q$-extension of exponential function as follows:

\begin{equation}\label{2}
\begin{split}
&e_q(z)=\sum_{n=0}^{\infty}\frac{z^n}{[n]_q!}=\frac{1}{((1-q)z:q)_{\infty}},~~(see  [3,6,7,8]),\\
\end{split}
\end{equation}

where $z \in \mathbb{C}~with~|z|<1$.

The Jackson definite $q$-integral of the function $f$ is defined by

\begin{equation}\label{3}
\int_{0}^{x}f(t)d_qt=(1-q)\sum_{a=0}^{\infty}f(q^ax)xq^a,~~see~[3,6,9]).
\end{equation}
 The $q$-defference operator $D_q$ is defined by

 \begin{equation}\label{4}
D_qf(x)=\frac{d_qf(x)}{d_qx}=
\begin{cases} \frac{f(x)-f(qx)}{(1-q)x} & \text{if $x \neq 0$}
\\
\frac{df(x)}{dx} & \text{if $x=0$,}
\end{cases}
\end{equation}

where

\begin{equation*}
\lim_{q \rightarrow 1}D_qf(x)=\frac{df(x)}{dx},~~\text{(see[3,6,8,10]).}
\end{equation*}

By using exponential function $e_q(x)$, Hegazi and Mausour defined $q$-Bernoulli polynomials by means of

\begin{equation}\label{5}
\sum_{n=0}^{\infty}B_{n,q}(x)\frac{t^n}{[n]_q!}=\frac{t}{e_q(t)-1}e_q(xt),~~\text{(see [3,6,8,12]).}
\end{equation}

In the special case, $x=0$, $B_{n,q}(0)=B_{n,q}$ are called the $n$-th $q$-Bernoulli numbers.  

From (\ref{5}), we can easily derive the following equation:

\begin{equation}\label{6}
B_{n,q}(x)=\sum_{l=0}^{n}\binom{n}{l}_qB_{n-l,q}x^l=\sum_{l=0}^{n}\binom{n}{l}_qB_{l,q}x^{n-l},
\end{equation}

where
\begin{equation*}
\binom{n}{l}_q=\frac{[n]_q!}{[n-l]_q![l]_q!}=\frac{[n]_q[n-1]_q \cdots [n-l+1]_q}{[l]_q!},~~\text{(see [6,12]).}
\end{equation*}

In the next section, we will consider new $q$-extensions of Euler numbers and polynomials by using the method of Hegazi and Mansour.  More than five decades ago,
Carlitz[2] defined a $q$-extension of Euler polynomials.  In a recent paper (see [7]), B. A. Kupershmidt constructed reflection symmetries of $q$-Bernoulli polynomials which
differ from Carlitz's $q$-Bernoulli numbers and polynomials.  By using the method of B. A. Kupershmidt, Hegazi and Mansour also introduced new $q$-extension of Bernoulli numbers and polynomials (see [3,7,8]).  From the $q$-exponential function, Kurt and Cenkci derived some interesting new formulae of $q$-extension of Genocchi polynomials.  Recently, several authors have studied various $q$-extension of Bernoulli and Euler polynomials (see [1-10]).  Let $\mathbb{C}$ be the complex number field and let $\mathcal{F}$ be the set of all formal power series in variable $t$ over $\mathbb{C}$ with

\begin{equation}\label{7}
\mathcal{F}=\left\{f(t)=\sum_{k=0}^{\infty}\frac{a_{k}}{[k]_q!}t^{k}\vert a_{k}\in\mathbb{C}\right\}.
\end{equation}

Let $\mathbb{P}=\mathbb{C}[t]$ and let $\mathbb{P}^*$ be the vector space of all linear
functionals on $\mathbb{P}$. $\langle L\vert p(x)\rangle$ denotes the action of linear functional $L$ on the polynomial $p(x)$, and it is well known that the vector
space operations on $\mathbb{P}^*$ are defined by

\begin{equation*}
\langle L+M \vert p(x)\rangle=\langle L\vert p(x)\rangle + \langle M\vert p(x)\rangle,~~\langle cL\vert p(x)\rangle=c\langle L\vert p(x)\rangle,
\end{equation*}

where $c$ is complex constant (see[1,5,12]).

For $f(t)=\sum_{k=0}^{\infty}\frac{a_{k}}{[k]_q!}t^{k} \in \mathcal{F}$, we define the linear functional on $\mathbb{P}$ by setting

\begin{equation}\label{8}
\langle f(t) \vert x^n \rangle = a_n~~\text{for all $n \geq 0$}.
\end{equation}

From (\ref{7}) and (\ref{8}), we note that

\begin{equation}\label{9}
\langle t^{k}\vert x^{n}\rangle =[n]_q!\delta_{n,k},\quad (n, k\geq 0),
\end{equation}

where $\delta_{n,k}$ is the Kronecker's symbol.\\

Let us assume that $f_{L}(t)=\sum_{k=0}^{\infty}\langle L\vert x^{n}\rangle\frac{t^{k}}{k!}$. Then by (\ref{9}), we easily see that
$\langle f_{L}(t)\vert x^{n}\rangle =\langle L\vert x^{n}\rangle$.
That is, $f_L(t)=L$. Additionally, the map $L\longmapsto f_{L}(t)$ is a vector space isomorphism from $\mathbb{P}^{*}$ onto $\mathcal{F}$. 
Henceforth $\mathcal{F}$ denotes both the algebra of formal power series in $t$ and the vector space of all linear functionals on $\mathbb{P}$, and so an element
$f(t)$ of $\mathcal{F}$ will be thought as a formal power series and a linear functional.  We call it the $q$-umbral algebra. The $q$-umbral calculus is the study of $q$-umbral
algebra. By (\ref{2}), and (\ref{3}),
we easily see that $\langle e_q(yt)\vert x^n\rangle =y^n$ and so $\langle e_q(yt)\vert p(x)\rangle =p(y)$ for $p(x) \in \mathbb{P}$.  The order $o(f(t))$ of the power series $f(t) \neq 0$ is the smallest integer for which $a_k$ does not vanish.  If $o(f(t))=0$, then $f(t)$ is called an invertible series.  If $o(f(t))=1$, then$f(t)$ is called a delta series (see [1,5,11,12]).  For $f(t),g(t) \in \mathcal{F}$, we have $\langle f(t)g(t)\vert p(x)\rangle =\langle f(t)\vert g(t)p(x)\rangle=\langle g(t)\vert f(t)p(x)\rangle$.
Let $f(t) \in \mathcal{F}$ and $p(x) \in \mathbb{P}$. Then we have

\begin{equation}\label{10}
f(t)=\sum_{k=0}^{\infty} \langle f(t)\vert x^{k}\rangle \frac{t^k}{[k]_q!},~~ p(x)=\sum_{k=0}^{\infty} \langle t^{k}\vert p(x)\rangle
\frac{x^k}{[k]_q!} ~~(\text{see [11]}).
\end{equation}

From (\ref{10}), we have

\begin{equation}\label{11}
p^{(k)}(x)=D^k_qp(x)=\sum_{l=k}^{\infty}
\frac{\langle t^{l}\vert p(x)\rangle}{[l]_q!}[l]_q \cdots [l-k+1]_qx^{l-k}.
\end{equation}

By (\ref{11}), we get

\begin{equation}\label{12}
p^{(k)}(0)=\langle t^{k}\vert p(x)\rangle ~~and~~\langle 1\vert p^{(k)}(x)\rangle=p^{(k)}(0).
\end{equation}

Thus from (\ref{12}), we note that

\begin{equation}\label{13}
t^{k}p(x)=p^{(k)}(x)=D^k_qp(x).
\end{equation}

Let $f(t),g(t) \in \mathcal{F}$ with $o(f(t))=1$ and $o(g(t))=0$. Then there exists a unique sequence $S_n(x)$(deg$S_n(x)=n$) of polynomials
such that $\langle g(t)f(t)^k \vert S_n(x)\rangle =[n]_q! \delta_{n,k},~~(n,k\geq 0)$.
The sequence $S_n(x)$ is called the $q$-sheffer sequence for $(g(t),f(t))$ which is denoted by $S_n(x) \sim (g(t),f(t))$.
Let $S_n(x) \sim (g(t),f(t))$. For $h(t) \in \mathcal{F}$ and $p(x) \in \mathbb{P}$, we have

\begin{equation}\label{14}
h(t)=\sum_{k=0}^{\infty}\frac{\langle h(t)\vert S_k(x)\rangle}{[k]_q!}g(t)f(t)^k,~~
p(x)=\sum_{k=0}^{\infty}
\frac{ \langle g(t)f(t)^k\vert p(x)\rangle}{[k]_q!}S_k(x),
\end{equation}

and

\begin{equation}\label{15}
\frac{1}{g(\bar{f}(t)}e_q(y \bar{f}(t))=\sum_{k=0}^{\infty}\frac{S_k(y)}{[k]_q!}t^k,~~\text{for all $y \in \mathbb{C}$},
\end{equation}

where $\bar{f}(t)$ is the compositional inverse of $f(t)$ (see [11,12]).

Recently, Araci-Acikgoz-Sen derived some new interesting properties on the new family of $q$-Euler numbers and polynomials from some applications of umbral algebra (see[1]).  The properties of $q$-Euler and 
$q$-Bernoulli polynomials seem to be of interest and worthwhile in the areas of both number theory and mathematical physics. In this paper, we develop the new method of $q$-umbral calculus due to Roman and study new $q$-extension of Euler numbers and polynomials which are derived from $q$-umbral calculus. Finally, we give new expicit formulas on $q$-Euler polynomials relate to Hegazi-Mansour's $q$-Bernoulli polynomials.

\section{$q$-Euler numbers and polynomials}
We consider the new $q$-extension of Euler polynomials which are generated by the generating function to be

\begin{equation}\label{16}
\frac{2}{e_q(t)+1}e_q(xt)=\sum_{n=0}^{\infty}E_{n,q}(x)\frac{t^n}{[n]_q!}.
\end{equation}

In the special case, $x=0$, $E_{n,q}(0)=E_{n,q}$ are called the $n$-th $q$-Euler numbers.
From (\ref{16}), we note that

\begin{equation}\label{17}
E_{n,q}(x)=\sum_{l=0}^{n}\binom{n}{l}_qE_{l,q}x^{n-l}=\sum_{l=0}^{n}\binom{n}{l}_qE_{n-l,q}x^{l}.
\end{equation}

By (\ref{16}), we easily get

\begin{equation}\label{18}
E_{0,q}=1,~~ E_{n,q}(1)+E_{n,q}=2 \delta_{0,n}
\end{equation}

For example, $E_{0,q}=1, E_{1,q}=-\frac{1}{2},~~E_{2,q}=\frac{q-1}{4},~~E_{3,q}=\frac{q+q^2-1}{4}+\frac{(1-q)[3]_q}{8},\cdots$.
From (\ref{15}) and (\ref{16}), we have

\begin{equation}\label{19}
E_{n,q}(x) \sim  \left(  \frac{e_q(t)+1}{2},t\right)
\end{equation}

and
\begin{equation}\label{20}
\frac{2}{e_q(t)+1}x^n=E_{n,q}(x),~~(n \geq 0).
\end{equation}

Thus, by (\ref{13}) and (\ref{20}), we get
\begin{equation}\label{21}
tE_{n,q}(x)=\frac{2}{e_q(t)+1}tx^n=[n]_q\frac{2}{e_q(t)+1}x^{n-1}=[n]_qE_{n-1,q}(x),~~(n \geq 0).
\end{equation}

Indeed, by (\ref{9}), we get

\begin{equation}\label{22}
\begin{split}
\langle \frac{e_q(t)+1}{2}t^k \vert E_{n,q}(x)\rangle= &\frac{[k]_q!}{2}\binom{n}{k}_q \langle e_q(t)+1 \vert E_{n-k,q}(x)\rangle\\
&=\frac{[k]_q!}{2}\binom{n}{k}_q ( E_{n-k,q}(1)+ E_{n-k,q}).\\
\end{split}
\end{equation}

From (\ref{19}), we have

\begin{equation}\label{23}
\langle \left( \frac{e_q(t)+1}{2}\right) t^k \vert E_{n,q}(x)\rangle= [n]_q! \delta_{n,k}
\end{equation}

Thus, by (\ref{22}) and (\ref{23}), we get

\begin{equation}\label{24}
0=E_{n-k,q}(1)+E_{n-k,q}=\sum_{l=0}^{n-k}\binom{n-k}{l}_qE_{l,q}+E_{n-k,q},~~(n,k \in \mathbb{Z}_{\geq 0}~ with~ n>k).
\end{equation}

This is equivalent to

\begin{equation}\label{25}
-2E_{n-k,q}=\sum_{l=0}^{n-k-1}\binom{n-k}{l}_qE_{l,q},~~where~~n,k \in \mathbb{Z}_{\geq 0}~ with ~n>k.
\end{equation}

Therefore , by (\ref{25}), we obtain the following lemma.

\begin{lemma}\label{lem1}
For $n \geq 1$, we have
\begin{equation*}
-2E_{n,q}=\sum_{l=0}^{n-1}\binom{n}{l}_qE_{l,q}.
\end{equation*}
\end{lemma}

From (\ref{17}) we have

\begin{equation}\label{26}
\begin{split}
\int _{x}^{x+y}E_{n,q}(u)d_qu&=\sum_{l=0}^{n}\binom{n}{l}_qE_{n-l,q} \frac{1}{[l+1]_q}\{(x+y)^{l+1}-x^{l+1}\}\\
&=\frac{1}{[n+1]_q}\sum_{l=0}^{n}\binom{n+1}{l+1}_qE_{n-l,q}\{(x+y)^{l+1}-x^{l+1}\}\\
&=\frac{1}{[n+1]_q}\sum_{l=1}^{n+1}\binom{n+1}{l}_qE_{n+1-l,q}\{(x+y)^{l}-x^{l}\}\\
&=\frac{1}{[n+1]_q}\sum_{l=0}^{n+1}\binom{n+1}{l}_qE_{n+1-l,q}\{(x+y)^{l}-x^{l}\}\\
&=\frac{1}{[n+1]_q}\{E_{n+1,q}(x+y)-E_{n+1,q}(x)\}.\\
\end{split}
\end{equation}

Thus, by (\ref{26}), we get
\begin{equation}\label{27}
\begin{split}
\langle \frac{e_q(t)-1}{t} \vert E_{n,q}(x)\rangle = &\frac{1}{[n+1]_q}\langle \frac{e_q(t)-1}{t} \vert tE_{n+1,q}(x)\rangle\\
= &\frac{1}{[n+1]_q}\langle e_q(t)-1 \vert E_{n+1,q}(x)\rangle\\
= &\frac{1}{[n+1]_q}\{ E_{n+1,q}(1)-E_{n+1,q} \}\\
=&\int _{0}^{1}E_{n,q}(u)d_qu.\\
\end{split}
\end{equation}

Therefore, by (\ref{27}), we obtain the following theorem.

\begin{theorem}\label{thm2}
For $n \geq 0$, we have
\begin{equation*}
\langle \frac{e_q(t)-1}{t} \vert E_{n,q}(x)\rangle =\int _{0}^{1}E_{n,q}(u)d_qu.
\end{equation*}
\end{theorem}

Let
\begin{equation}\label{28}
\mathbb{P}_n=\{p(x) \in \mathbb{C}[x] \vert deg~ p(x) \leq n \}.
\end{equation}

For $p(x) \in \mathbb{P}_n$, let us assume  that

\begin{equation}\label{29}
p(x)=\sum_{k=0}^{n}b_{k,q}E_{k,q}(x).
\end{equation}

Then, by (\ref{19}), we get

\begin{equation}\label{30}
\langle \left( \frac{e_q(t)+1}{2}\right) t^k \vert E_{n,q}(x)\rangle= [n]_q! \delta_{n,k}.
\end{equation}

From (\ref{29}) and (\ref{30}), we can derive the following equation (\ref{31}):
\begin{equation}\label{31}
\begin{split}
\langle \left( \frac{e_q(t)+1}{2}\right) t^k \vert p(x)\rangle &= \sum_{l=0}^{n}b_{l,q}
\langle \left( \frac{e_q(t)+1}{2}\right) t^k \vert E_{l,q}(x)\rangle\\
&=\sum_{l=0}^{n}b_{l,q}[l]_q!\delta_{l,k}=[k]_q!b_{k,q}.\\
\end{split}
\end{equation}

Thus, by (\ref{31}), we get

\begin{equation}\label{32}
\begin{split}
b_{k,q}&=\frac{1}{[k]_q!}\langle \left( \frac{e_q(t)+1}{2}\right) t^k \vert p(x)\rangle
=\frac{1}{2[k]_q!}\langle \left( e_q(t)+1 \right) t^k \vert p(x)\rangle\\
&=\frac{1}{2[k]_q!}\langle e_q(t)+1 \vert p^{(k)}(x)\rangle
=\frac{1}{2[k]_q!} \{ p^{(k)}(1)+ p^{(k)}(0)\},\\
\end{split}
\end{equation}

where $p^{(k)}(x)=D_q^{k}p(x)$.

Therefore, by (\ref{29}) and (\ref{32}), we obtain the following theorem.
\begin{theorem}\label{thm3}
For $p(x) \in \mathbb{P}_n$, let $p(x)=\sum_{k=0}^{n}b_{k,q}E_{k,q}(x).$
Then we have\\
\begin{equation*}
\begin{split}
b_{k,q}&=\frac{1}{2[k]_q!}\langle \left( e_q(t)+1 \right) t^k \vert p(x)\rangle\\
&=\frac{1}{2[k]_q!} \{ p^{(k)}(1)+ p^{(k)}(0)\},\\
\end{split}
\end{equation*}
\end{theorem}

where $p^{(k)}(x)=D_q^kp(x)$.\\

From (\ref{5}), we note that
\begin{equation}\label{33}
B_{n,q}(x) \sim  \left(  \frac{e_q(t)-1}{t},t\right),~~(n \geq 0).
\end{equation}

Let us take $p(x)=B_{n,q}(x)\in \mathbb{P}_n$.  Then $B_{n,q}(x)$ can be represented as a linear combination of $\{E_{0,q}(x),E_{1,q}(x)
, \cdots , E_{n,q}(x)\}$ as follows:

\begin{equation}\label{34}
B_{n,q}(x)=p(x)=\sum_{k=0}^{n}b_{k,q}E_{k,q}(x),~~(n \geq 0),
\end{equation}

where
\begin{equation}\label{35}
\begin{split}
b_{k,q}&=\frac{1}{2[k]_q!}\langle \left( e_q(t)+1\right) t^k \vert B_{n,q}(x)\rangle\\
&=\frac{[n]_q[n-1]_q \cdots [n-k+1]_q}{2[k]_q!}\langle  e_q(t)+1 \vert B_{n-k,q}(x)\rangle\\
&=\frac{1}{2} \binom{n}{k}_q\langle  e_q(t)+1  \vert B_{n-k,q}(x)\rangle
=\frac{1}{2} \binom{n}{k}_q \{ B_{n-k,q}(1)+B_{n-k,q} \}.\\
\end{split}
\end{equation}

From (\ref{5}), we can derive the following recurrence relation for the $q$-Bernoulli numbers:

\begin{equation}\label{36}
\begin{split}
t&=\left( \sum_{l=0}^{\infty}B_{l,q}\frac{t^l}{[l]_q!} \right)(e_q(t)-1)\\
&=\sum_{n=0}^{\infty}\left( \sum_{l=0}^{n}\binom{n}{l}_q
B_{l,q} \right)\frac{t^n}{[n]_q!}-\sum_{n=0}^{\infty}B_{n,q}\frac{t^n}{[n]_q!}\\
&=\sum_{n=0}^{\infty}\left(B_{n,q}(1)-B_{n,q} \right)\frac{t^n}{[n]_q!}.\\
\end{split}
\end{equation}

Thus, by (\ref{36}), we get

\begin{equation}\label{37}
B_{0,q}=1, B_{n,q}(1)-B_{n,q}=
\begin{cases} 1 & \text{if $n=1$,}
\\
0 & \text{if $n>1$.}
\end{cases}
\end{equation}

For example, $B_{0,q}=1,B_{1,q}=-\frac{1}{[2]_q},~~B_{2,q}=\frac{q^2}{[3]_q[2]_q},\cdots$.\\

By (\ref{34}), (\ref{35}), and (\ref{37}), we get

\begin{equation}\label{38}
\begin{split}
B_{n,q}(x)&=b_{n,q}E_{n,q}(x)+b_{n-1,q}E_{n-1,q}(x)+\sum_{k=0}^{n-2}b_{k,q}E_{k,q}(x)\\
&=E_{n,q}(x)+\frac{[n]_q}{2}\left( 1-\frac{2}{[2]_q} \right)E_{n-1,q}(x)+\sum_{k=0}^{n-2}\binom{n}{k}_qB_{n-k,q}E_{k,q}(x)\\
&=E_{n,q}(x)-\frac{[n]_q(1-q)}{2[2]_q}E_{n-1,q}(x)+\sum_{k=0}^{n-2}\binom{n}{k}_qB_{n-k,q}E_{k,q}(x)\\
\end{split}
\end{equation}

Therefore, by (\ref{38}), we obtain the following theorem.
\begin{theorem}\label{thm4}
For $n \geq 2$, we have
\begin{equation*}
B_{n,q}(x)=E_{n,q}(x)+\frac{[n]_q(q-1)}{2[2]_q}E_{n-1,q}(x)+\sum_{k=0}^{n-2}\binom{n}{k}_qB_{n-k,q}E_{k,q}(x).
\end{equation*}
\end{theorem}

For $r \in \mathbb{Z}_{\geq 0}$, the $q$-Euler polynomials, $E_{n,q}^{(r)}(x)$, of order $r$ are defined by the generating function to be

\begin{equation}\label{39}
\begin{split}
\left(  \frac{2}{e_q(t)+1}\right)^r e_q(xt)&=\underbrace{\left(  \frac{2}{e_q(t)+1}\right) \times  \cdots \times \left(  \frac{2}{e_q(t)+1}\right)}_{r-times}e_q(xt)\\
&=\sum_{n=0}^{\infty}E_{n,q}^{(r)}(x)\frac{t^n}{[n]_q!}.\\
\end{split}
\end{equation}

In the special case, $x=0,E_{n,q}^{(r)}(0)=E_{n,q}^{(r)} $ are called the $n$-th $q$-Euler numbers of order $r$.\\

Let
\begin{equation}\label{40}
g^r(t)=\left(  \frac{e_q(t)+1}{2}\right)^r, ~~(r \in \mathbb{Z}_{\geq 0}).
\end{equation}

Then $g^r(t)$ is an invertible series.
From (\ref{39}) and (\ref{40}), we have

\begin{equation}\label{41}
\sum_{n=0}^{\infty}E_{n,q}^{(r)}(x)\frac{t^n}{[n]_q!}=\frac{1}{g^r(t)}e_q(xt)=\sum_{n=0}^{\infty}\frac{1}{g^r(t)}x^n\frac{t^n}{[n]_q!}.
\end{equation}

By (\ref{41}), we get
\begin{equation}\label{42}
E_{n,q}^{(r)}(x)=\frac{1}{g^r(t)}x^n,
\end{equation}

and

\begin{equation}\label{43}
tE_{n,q}^{(r)}(x)=\frac{1}{g^r(t)}tx^n=[n]_q \frac{1}{g^r(t)}x^{n-1}=[n]_qE_{n-1,q}^{(r)}(x).
\end{equation}

Thus, by (\ref{41}), (\ref{42}) and (\ref{43}), we see that

\begin{equation}\label{44}
E_{n,q}^{(r)}(x) \sim  \left(  \left(\frac{e_q(t)+1}{2}\right)^r,t\right).
\end{equation}

By (\ref{9}) and (\ref{39}), we get
\begin{equation}\label{45}
\langle \left(\frac{2}{e_q(t)+1}\right)^re_q(yt) \vert x^n \rangle= E_{n,q}^{(r)}(y)
=\sum_{l=0}^{n}\binom{n}{l}_q E_{n-l,q}^{(r)}y^l.
\end{equation}

Thus, we have
\begin{equation}\label{46}
\begin{split}
\langle \left(\frac{2}{e_q(t)+1}\right)^r \vert x^n \rangle&=\sum_{m=0}^{\infty}\left( \sum_{i_1+ \cdots +i_r=m}^{}
\frac{E_{i_1,q} \cdots E_{i_r,q}}{[i_1]_q! \cdots [i_r]_q!} \right) \langle t^m \vert x^n \rangle\\
&=\sum_{i_1+ \cdots +i_r=n}\frac{[n]_q!}{[i_1]_q! \cdots [i_r]_q!}E_{i_1,q} \cdots E_{i_r,q}\\
&=\sum_{i_1+ \cdots +i_r=n}\binom{n}{i_1, \cdots , i_r}_q E_{i_1,q} \cdots E_{i_r,q},\\
\end{split}
\end{equation}

where $\binom{n}{i_1, \cdots , i_r}_q=\frac{[n]_q!}{[i_1]_q! \cdots [i_r]_q!}$.\\

By (\ref{45}), we easily get
\begin{equation}\label{47}
\langle \left(\frac{2}{e_q(t)+1}\right)^r \vert x^n \rangle= E_{n,q}^{(r)}.
\end{equation}

Therefore, by (\ref{46}) and (\ref{47}), we obtain the following theorem.
\begin{theorem}\label{thm5}
For $n \geq 0$, we have
\begin{equation*}
E_{n,q}^{(r)}=\sum_{i_1+ \cdots +i_r=n}\binom{n}{i_1, \cdots , i_r}_q E_{i_1,q} \cdots E_{i_r,q},
\end{equation*}
\end{theorem}

where $\binom{n}{i_1, \cdots , i_r}_q=\frac{[n]_q!}{[i_1]_q! \cdots [i_r]_q!}$.

Let us take $p(x)=E_{n,q}^{(r)}(x) \in \mathbb{P}_n$. Then, by Theorem \ref{thm3}, we get

\begin{equation}\label{48}
E_{n,q}^{(r)}(x)=p(x)=\sum_{k=0}^{n}b_{k,q}E_{k,q}(x),
\end{equation}

where
\begin{equation}\label{49}
\begin{split}
b_{k,q}&=\frac{1}{2[k]_q!}\langle \left( e_q(t)+1\right) t^k \vert p(x)\rangle=\frac{1}{2[k]_q!}\langle \left( e_q(t)+1\right) \vert t^k p(x)\rangle\\
&=\frac{\binom{n}{k}_q}{2}\langle \left( e_q(t)+1\right)\vert E_{n-k,q}^{(r)}(x)\rangle=\frac{\binom{n}{k}_q}{2}\{E_{n-k,q}^{(r)}(1)+E_{n-k,q}^{(r)}\}.
\end{split}
\end{equation}

From (\ref{39}), we have

\begin{equation}\label{50}
\begin{split}
\sum_{k=0}^{\infty}& \{E_{n,q}^{(r)}(1)+E_{n,q}^{(r)}\}\frac{t^n}{[n]_q!}=\left(\frac{2}{e_q(t)+1}\right)^r(e_q(t)+1)\\
&=2 \left(\frac{2}{e_q(t)+1}\right)^{r-1}=2\sum_{n=0}^{\infty}E_{n,q}^{(r-1)}\frac{t^n}{[n]_q!}.
\end{split}
\end{equation}

By comparing the coefficients on the both sides of (\ref{50}), we get
\begin{equation}\label{51}
E_{n,q}^{(r)}(1)+E_{n,q}^{(r)}=2E_{n,q}^{(r-1)}, ~~(n \geq 0).
\end{equation}

Therefore, by (\ref{48}), (\ref{49}) and (\ref{51}), we obtain the following theorem.
\begin{theorem}\label{thm6}
For $n\in \mathbb{Z}_{\geq 0}$,  $r\in \mathbb{Z}_{> 0}$, we have
\begin{equation*}
E_{n,q}^{(r)}(x)=\sum_{k=0}^{\infty}\binom{n}{k}_q E_{n-k,q}^{(r-1)}E_{k,q}(x).
\end{equation*}
\end{theorem}

Let us assume that

\begin{equation}\label{52}
p(x)=\sum_{k=0}^{n}b_{k,q}^{r}E_{k,q}^{(r)}(x)\in \mathbb{P}_n.
\end{equation}

By (\ref{44}) and (\ref{52}), we get

\begin{equation}\label{53}
\begin{split}
\langle \left(\frac{e_q(t)+1}{2}\right)^r t^k \vert p(x) \rangle &=\sum_{l=0}^{n}b_{l,q}^{r}\langle \left(\frac{e_q(t)+1}{2}\right)^r t^k
\vert E_{l,q}^{(r)}(x) \rangle\\
&=\sum_{l=0}^{n}b_{l,q}^{r}[l]_q!\delta_{l,k}=[k]_q!b_{k,q}^{r}.\\
\end{split}
\end{equation}

From (\ref{53}), we have

\begin{equation}\label{54}
\begin{split}
b_{k,q}^{r}&=\frac{1}{[k]_q!}\langle \left(\frac{e_q(t)+1}{2}\right)^r t^k \vert p(x) \rangle=\frac{1}{2^r[k]_q!}
\langle \left( e_q(t)+1\right)^r \vert t^k p(x)\rangle\\
&=\frac{1}{2^r[k]_q!}\sum_{l=0}^{r}\binom{r}{l}\sum_{m \geq 0} \left( \sum_{i_1+ \cdots +i_l=m}\binom{m}{i_1, \cdots , i_l}_q\right)
\frac{1}{[m]_q!}\langle 1 \vert t^{m+k} p(x) \rangle\\
&=\frac{1}{2^r[k]_q!}\sum_{l=0}^{r}\binom{r}{l}\sum_{m \geq 0}  \sum_{i_1+ \cdots +i_l=m}\binom{m}{i_1, \cdots , i_l}_q\frac{1}{[m]_q!}
p^{(m+k)}(0).
\end{split}
\end{equation}

Therefore by (\ref{52}) and (\ref{54}), we obtain the following theorem.

\begin{theorem}\label{thm7}
For $n \geq 0$, let
$p(x)=\sum_{k=0}^{n}b_{k,q}^{r}E_{k,q}^{(r)}(x)\in \mathbb{P}_n$.\\
Then we have\\
\begin{equation*}
\begin{split}
b_{k,q}^{r}&=\frac{1}{2^r[k]_q!}=\langle \left( e_q(t)+1\right)^r t^k \vert p(x)\rangle\\
&=\frac{1}{2^r[k]_q!}\sum_{m \geq 0}\sum_{l=0}^{r}\binom{r}{l}  \sum_{i_1+ \cdots +i_l=m}\binom{m}{i_1, \cdots , i_l}_q \frac{1}{[m]_q!}
p^{(m+k)}(0),\\
\end{split}
\end{equation*}
\end{theorem}

where $p^{(k)}(x)=D^k_q p(x)$.

Let us take $p(x)=E_{n,q}(x)\in \mathbb{P}_n$. Then, by Theorem \ref{thm7}, we get
\begin{equation}\label{55}
E_{n,q}(x)=p(x)=\sum_{k=0}^{n}b_{k,q}^{r}E_{k,q}^{(r)}(x),
\end{equation}

where

\begin{equation}\label{56}
\begin{split}
b_{k,q}&=\frac{1}{2^r[k]_q!}\sum_{m=0}^{n-k}\sum_{l=0}^{r}\binom{r}{l} \sum_{i_1+ \cdots +i_l=m}\binom{m}{i_1, \cdots , i_l}_q\\
&\times\;\frac{1}{[m]_q!}[n]_q \cdots [n-m-k+1]_qE_{n-m-k,q}\\
&=\frac{1}{2^r}\sum_{m=0}^{n-k}\sum_{l=0}^{r}\binom{r}{l} \sum_{i_1+ \cdots +i_l=m}\binom{m}{i_1, \cdots , i_l}_q\\
&\times\;\frac{[m+k]_q!}{[m]_q![k]_q!}\frac{[n]_q \cdots [n-m-k+1]_q}{[m+k]_q!}E_{n-m-k,q}\\
&=\frac{1}{2^r}\sum_{m=0}^{n-k}\sum_{l=0}^{r} \sum_{i_1+ \cdots +i_l=m}\binom{r}{l}\binom{m}{i_1, \cdots , i_l}_q \binom{m+k}{m}_q \binom{n}{m+k}_q E_{n-m-k,q}\\
\end{split}
\end{equation}

Therefore, by (\ref{55}) and (\ref{56}), we obtain the following theorem.

\begin{theorem}\label{thm8}
For $n,r \geq 0$, we have
\begin{eqnarray*}
E_{n,q}(x)=\frac{1}{2^r}\sum_{k=0}^{n} \Bigg\{ \sum_{m=0}^{n-k}\sum_{l=0}^{r} \sum_{i_1+ \cdots +i_l=m}\binom{r}{l}\binom{m}{i_1, \cdots , i_l}_q \binom{m+k}{m}_q
\binom{n}{m+k}_q
\end{eqnarray*}
\begin{eqnarray*}
\times\;E_{n-m-k,q} \Bigg\} E_{k,q}^{(r)}(x).
\end{eqnarray*}
\end{theorem}
For $r \in \mathbb{Z}_{\geq 0}$, let us consider $q$-Bernoulli polynomials of order $r$ which are defined by the generating function to be

\begin{equation}\label{57}
\begin{split}
\left(  \frac{t}{e_q(t)-1}\right)^r e_q(xt)&=\underbrace{\left(  \frac{t}{e_q(t)-1}\right) \times  \cdots \times \left(  \frac{t}{e_q(t)-1}\right)}_{r-times}e_q(xt)\\
&=\sum_{n=0}^{\infty}B_{n,q}^{(r)}(x)\frac{t^n}{[n]_q!}.\\
\end{split}
\end{equation}

In the special case, $x=0$, $B_{n,q}^{(r)}(0)=B_{n,q}^{(r)}$ are called the $n$-th $q$-Bernoulli numbers of order $r$.
By (\ref{57}), we easily get

\begin{equation}\label{58}
B_{n,q}^{(r)}(x)=\sum_{l=0}^{n}\binom{n}{l}_q B_{l,q}^{(r)}x^{n-l} \in \mathbb{P}_n.
\end{equation}

Let us take $p(x)=B_{n,q}^{(r)}(x)\in \mathbb{P}_n$.  Then, by Theorem \ref{thm7}, we get

\begin{equation}\label{59}
B_{n,q}^{(r)}(x)=p(x)=\sum_{k=0}^{n}b_{k,q}^{r}E_{k,q}^{(r)}(x),
\end{equation}
where

\begin{equation}\label{60}
\begin{split}
b_{k,q}^{r}&=\frac{1}{2^r[k]_q!}\langle \left( e_q(t)+1\right)^r t^k \vert B_{n,q}^{(r)}(x))\rangle\\
&=\frac{1}{2^r[k]_q!}\sum_{m=0}^{n-k}\sum_{l=0}^{r}\binom{r}{l} \sum_{i_1+ \cdots +i_l=m}\binom{m}{i_1, \cdots , i_l}_q\\
&\times\;\frac{[n]_q \cdots [n-m-k+1]_q}{[m]_q!}B_{n-m-k,q}^{(r)}\\
&=\frac{1}{2^r}\sum_{m=0}^{n-k}\sum_{l=0}^{r} \sum_{i_1+ \cdots +i_l=m}\binom{r}{l}\binom{m}{i_1, \cdots , i_l}_q \binom{m+k}{m}_q \binom{n}{m+k}_q
B_{n-m-k,q}^{(r)}.\\
\end{split}
\end{equation}

Therefore, by (\ref{59}) and (\ref{60}), we obtain the following theorem.

\begin{theorem}\label{thm9}
For $n,r \geq 0$, we have
\begin{eqnarray*}
B_{n,q}^{(r)}(x)=\frac{1}{2^r}\sum_{k=0}^{n} \Bigg\{ \sum_{m=0}^{n-k}\sum_{l=0}^{r} \sum_{i_1+ \cdots +i_l=m}\binom{r}{l}\binom{m}{i_1, \cdots , i_l}_q
\binom{m+k}{m}_q \binom{n}{m+k}_q
\end{eqnarray*}
\begin{eqnarray*}
\times\;B_{n-m-k,q}^{(r)} \Bigg\} E_{k,q}^{(r)}(x).
\end{eqnarray*}
\end{theorem}

%%%%%%%%%%%%%%%%%%%%%%%%%%%%%%%%%%%%%%%%%%%%
%%%%%%%%%%%%%%%%%%%%%%%%%%%%%%%%%%%%%%%%%%%%
%%%%%%%%%%%%%%%%%%%%%%%%%%%%%%%%%%%%%%%%%%%%
%%%%%%%%%%%%%%%%%%%%%%%%%%%%%%%%%%%%%%%%%%%%
%%%%%%%%%%%%%%%%%%%%%%%%%%%%%%%%%%%%%%%%%%%%

\section*{Acknowledgements}

This work was supported by the National Research Foundation of Korea(NRF) grant funded by the Korea government(MOE) (No.2012R1A1A2003786 ).

\par\bigskip

\par

\bigskip\bigskip

\begin{center}\begin{large}

\end{large}\end{center}

\par

%%%%%%%%%%%%%%%%%%%%%%%%%%%%%%%%%%%%%%%%%%%%%%%%%%%%%%%%%%%%%%%%%%%%%%%%%%%%%%%%%%%%%%%%%%%%%%%%%%%%%%%%%%%%%%%%%%%%%%%%%%%%%%%%%%%%%%%%%%%%%%%%%%%%%%\begin{thebibliography}{99}

\end{document}